\documentclass[12pt]{article}
\usepackage[a4paper]{geometry}
\usepackage{amsthm}
\usepackage{amsmath}
\usepackage{amssymb}
\usepackage{amsfonts}

\def\N{\mathbb{N}}

\def\R{\mathbb{R}}

\def\E{\mathbb{E}}
\def\T{\mathbb{T}}

\def\Z{\mathbb{Z}}

\def\<{\big\langle}
\def\>{\big\rangle}
\def\Osc{\operatorname{Osc}}

\def\Var{\operatorname{Var}}
\def\Cov{\operatorname{Cov}}

\newtheorem{lemma}{Lemma}[section]

\newtheorem{theorem}{Theorem}[section]
\newtheorem{proposition}{Proposition}[section]

\theoremstyle{remark}

\theoremstyle{definition}

\newtheorem{definition}{Definition}[section]

\def\eref#1{(\ref{#1})}

\newtheorem{Hypothesis}{Hypothesis}

\begin{document}
\title{Super-diffusivity in a shear flow model from perpetual homogenization.}         
\author{G\'{e}rard Ben Arous\footnote{EPFL, gerard.benarous@epfl.ch} and Houman Owhadi\footnote{
Technion, owhadi@techunix.technion.ac.il}}     
\date{\today}          
\maketitle
\begin{abstract}
This paper is concerned with the asymptotic behavior solutions of  stochastic differential equations
$dy_t=d\omega_t -\nabla \Gamma(y_t) dt$, $y_0=0$ and $d=2$.  $\Gamma$ is a $2\times 2$ skew-symmetric matrix associated to a shear flow  characterized by an infinite number
 of spatial scales $\Gamma_{12}=-\Gamma_{21}=h(x_1)$, with $h(x_1)=\sum_{n=0}^\infty \gamma_n h^n(x_1/R_n)$ where $h^n$ are smooth functions of period 1, $h^n(0)=0$,  $\gamma_n$ and
$R_n$ grow exponentially fast with $n$. We can show that $y_t$ has an anomalous fast behavior ($\E[|y_t|^2]\sim t^{1+\nu}$ with $\nu>0$)
 and obtain quantitative estimates on the anomaly using and developing the tools of homogenization.
\end{abstract}

\tableofcontents

\section{Introduction}
Turbulent incompressible flows are characterized by multiple scales of mixing length and convection rolls. It is heuristically known and expected  that a diffusive transport in such media will be super-diffusive.
The first known observation of this anomaly is attributed to Richardson \cite{Ric26} who analyzed available experimental data on diffusion in air, varying  on about $12$ orders of magnitude. On that basis, he empirically conjectured that the diffusion coefficient $D_\lambda$ in turbulent air depends on the scale length $\lambda$ of the measurement. The Richardson law,
\begin{equation}
D_\lambda \propto \lambda^\frac{4}{3}
\end{equation}
was related to Kolmogorov-Obukhov turbulence spectrum, $v \propto \lambda^\frac{1}{3}$, by Batchelor \cite{Bat52}. The super-diffusive law of the root-mean-square relative displacement $\lambda(t)$ of advected particles
\begin{equation}
\lambda(t) \propto (D_{\lambda(t)}t)^\frac{1}{2}\propto t^\frac{3}{2}
\end{equation}
was derived by Obukhov \cite{Obu41} from a dimensional analysis similar to the one that led Kolmogorov \cite{Kol41}
to the $\lambda^\frac{1}{3}$ velocity spectrum.\\
More recently physicists and mathematicians have started to investigate on the super-diffusive phenomenon (from both heuristic and rigorous point of view) by using the tools of homogenization or renormalization; we refer to  M. Avellaneda and A. Majda \cite{AvMa90}, \cite{AvMa94}, J. Glimm and Al. \cite{FuGl90}, \cite{FuGl91}, \cite{Gli92}, J. Glimm and Q. Zhang \cite{GlZh92}, Q. Zhang \cite{Zh92}, M.B. Isichenko and J. Kalda \cite{IsKa91}, G. Gaudron \cite{Gau98}.

\nocite{Owh01}
It is now well known that homogenization over a periodic or ergodic divergence free drift has the property to enhance the diffusion \cite{FaPa94},\cite{FaPa96}, \cite{KoOl00}. It is also expected that several spatial scales of eddies should give rise to anomalous diffusion between proper time scales outside homogenization regime or when the bigger scale has not yet been homogenized. We refer  to
  M. Avellaneda \cite{Ave96}; A. Fannjiang \cite{Fan99}; Rabi Bhattacharya  \cite{Bha99} (see also \cite{BhDeGo99} by  Bhattacharya - Denker and  Goswami); A. Fannjiang and T. Komorowski \cite{FaKo01} and this panorama is certainly not complete.\\
The purpose of this paper is to implement rigorously on a shear flow model the idea that the key to anomalous fast diffusion in turbulent flows is an unfinished homogenization process over a large number of scales of eddies without a sharp separation between them. We will assume that the ratios between the spacial scales are bounded. The underlying phenomenon is similar to the one related to anomalous slow diffusion from perpetual homogenization on an infinite number of scales of gradient drifts \cite{Ow00a}, \cite{BeOw00b}, the main difference lies in the asymptotic behavior of the multi-scale effective diffusivities $D(n)$ associated with $n$ spacial scales, i.e. $D(n)$ diverge towards $\infty$ or converge towards $0$ with exponential rate depending on the nature of the scales: eddies or obstacles.\\
Note that the shear-layer model is exactly
solvable (\cite{AvMa90}, \cite{Gau98}).
When the geometrically divergent scales are recast into the Fourier
setting with a power-law spectrum, super-diffusivity has already been proven in the limit $t\rightarrow \infty$.
Our purpose in this paper is to show that never-ending homogenization can be used as a tool to obtain a quantitative control on the anomaly for finite times, not just an asymptotic result and without any self-similarity assumption. Moreover it will be shown that the mean-squared displacement $\E[y_t^2]$ of the diffusion in the shear flow behaves like $D(n(t)) t$ (see \eref{jhsasiuiuw798761}). In this formula,  $n(t)$ has a logarithmic growth and corresponds to the number of scales that can be considered as homogenized at time $t$, casting into light the role of never-ending homogenization in the anomalous fast behavior of a diffusion process in a shear-flow model. Moreover it will be shown in \cite{Ow00b} that the strategy associated to never-ending homogenization can be extended to higher dimensions (and non shear flow models of turbulence). Note added in proof: we would like to refer the reader to an interesting and related recent preprint by S. Olla and T. Komorowski \cite{OlKo} on "the Superdiffusive Behavior of Passive Tracer with a Gaussian Drift".

\section{The model}
Let us consider in dimension two a Brownian motion with a drift given by the divergence of a shear flow stream matrix, i.e. the solution of the stochastic differential equation:
\begin{equation}\label{SFsuperdiffstochdiffequ}
  dy_t = d\omega_t - \nabla \Gamma (y_t) dt , \quad y_0=0
\end{equation}
where $\Gamma$ is a skew-symmetric $2 \times 2$ shear flow matrix.
\begin{equation}
\Gamma(x_1,x_2)=\begin{pmatrix}
0 & h(x_1) \\
-h(x_1) & 0
\end{pmatrix}
\end{equation}
The function $(x_1,x_2)\rightarrow h(x_1)$ is given by a sum of infinitely many periodic functions with (geometrically) increasing periods
\begin{equation}\label{SFsupfracgraminfty}
h(x_1)=\sum_{n=0}^\infty \gamma_n h_n(\frac{x_1}{R_n})
\end{equation}
Where $h_n$ are smooth functions of period $1$. We will assume that
\begin{equation}\label{askdjjjj18}
h_n(0)=0
\end{equation}
We will normalize  the functions $h_n$ by the choosing their variance equal to one:
\begin{equation}\label{SFsupallngamnzeroo}
\Var(h_n)=\int_0^1 (h_n(x)-\int_0^1 h_n(y)dy)^2 dx=1
\end{equation}
$R_n$ and $\gamma_n$ grow exponentially fast with $n$, i.e.
\begin{equation}
R_n=\prod_{k=0}^n r_k
\end{equation}
Where $r_n$ are integers, $r_0=1$,
\begin{equation}\label{Modsubboundrnrhonmin}
 \rho_{\min}=\inf_{n\in \N^*} r_n \geq 2\quad\text{and}\quad \rho_{\max}=\sup_{n\in \N^*} r_n  < \infty
\end{equation}
We choose $\gamma_0=1$ and
\begin{equation}\label{Modsubboundrnrhonmin002}
 \gamma_{\min}=\inf_{n\in \N} (\gamma_{n+1}/\gamma_n) >1 \quad\text{and}\quad \gamma_{\max}=\sup_{n\in \N} (\gamma_{n+1}/\gamma_n)  < \infty
\end{equation}
It is assumed that  the first derivate of the potentials $h_n$ are uniformly bounded. ($\Osc(h)$ stands for $\sup h - \inf h$)
\begin{equation}\label{ModsubContUngradUn}
K_0=\sup_{n\in \N} \Osc(h_n) <\infty ,\quad \quad K_1=\sup_{n\in \N} \|h_n'\|_\infty <\infty
\end{equation}
In this paper we shall distinguish two hypotheses
\begin{Hypothesis}\label{csdcssjhxnjj9jdxc7h}
\begin{equation}\label{akjsbdu8191}
\rho_{\min}>\gamma_{\max}
\end{equation}
\end{Hypothesis}

\begin{Hypothesis}\label{csdcssjhxnjj9jdxc7h2}
\begin{equation}\label{akjsbdu8191lkkl2}
\rho_{\min}>\gamma_{\max}^{1/2}
\end{equation}
For all $n\in \N$
\begin{equation}\label{akjsbdu8191lkkl22s1}
h_n'(0)=0
\end{equation}
and
\begin{equation}\label{ModsubContUngradUndxhh2}
 K_2=\sup_{n\in \N} \|\partial_1^2 h_n\|_\infty <\infty
\end{equation}
\end{Hypothesis}

Let us observe that under the hypotheses \ref{csdcssjhxnjj9jdxc7h}, \eref{askdjjjj18}, \eref{Modsubboundrnrhonmin}, \eref{Modsubboundrnrhonmin002} and \eref{ModsubContUngradUn}   $h$ is a well defined $C^1$ function on $\R^d$ and
\begin{equation}
|h(x)|\leq K_1 |x| (1-\gamma_{\max}/\rho_{\min})^{-1}\quad |h'(x)|\leq K_1  (1-\gamma_{\max}/\rho_{\min})^{-1}
\end{equation}
Thus under the hypothesis \ref{csdcssjhxnjj9jdxc7h},
 $\Gamma$ is a well defined Lipschitz stream matrix  and the solution of the stochastic differential equation \eref{SFsuperdiffstochdiffequ} exists; is unique up to sets of measure $0$ with respect to the Wiener measure and is a strong Markov continuous Feller process.\\
Under the hypotheses \ref{csdcssjhxnjj9jdxc7h2}, \eref{askdjjjj18}, \eref{Modsubboundrnrhonmin}, \eref{Modsubboundrnrhonmin002} and \eref{ModsubContUngradUn}, $\Gamma$ is no more Lipschitz continuous but    $h$ is still a well defined $C^2$ function on $\R^d$ and
\begin{equation}\label{skjhxcloiw}
|h(x)|\leq K_2 |x|^2 (1-\gamma_{\max}/\rho_{\min}^2)^{-1}\quad |h'(x)|\leq K_2  (1-\gamma_{\max}/\rho_{\min}^2)^{-1} |x|
\end{equation}

\section{Main results}
\subsection{Under the hypothesis \ref{csdcssjhxnjj9jdxc7h}}
Our objective is to show that the solution \eref{SFsuperdiffstochdiffequ} is abnormally fast and the asymptotic sub-diffusivity will be characterized as an anomalous behavior of the variance at time $t$,
  i.e.  $\E_0[y_t^2] \sim t^{1+\nu}$ as $t\to\infty$.
More precisely  there exists a constant $\rho_0(\gamma_{min},\gamma_{max},K_0,K_1)$ and a time $t_0(\gamma_{min},\gamma_{max},R_1,K_0,K_1)$ such that
\begin{theorem}\label{IntSUMSFDThmsqdishgc1}
If $\rho_{min}>\rho_0$ and $y_t$ is a solution of \eref{SFsuperdiffstochdiffequ} then for $t>t_0$
\begin{equation}\label{hgdvsdckijeq71}
\E_0[|y_t.e_2|^2]=t^{1+\nu(t)}
\end{equation}
with
\begin{equation}
\frac{\ln \gamma_{\min}}{\ln \rho_{max}+\ln \frac{\gamma_{\max}}{\gamma_{\min}}}-\frac{C_1}{\ln t} \leq \nu(t)\leq \frac{\ln \gamma_{\max}}{\ln \rho_{min}+\ln \frac{\gamma_{\min}}{\gamma_{\max}}}+\frac{C_2}{\ln t}
\end{equation}
Where the constants $C_1$ and $C_2$ depends on $\rho_{\min},\gamma_{\min},\gamma_{\max},\rho_{\max},K_0,K_1$
\end{theorem}
We remark that if $\gamma_{\max}=\gamma_{\min}=\gamma$ and $\rho_{\max}=\rho_{\min}=\rho$ then $\nu(t)\sim \ln \gamma /\ln \rho$.\\
The key of the fast asymptotic behavior of the variance of the solution of \eref{SFsuperdiffstochdiffequ} is the geometric rate of divergence towards $\infty$ of the multi-scale effective matrices associated to a finite number of scales. More precisely, for $k,p\in \N$, $k\leq p$ we will write
\begin{equation}\label{hvjshvxvajkko1}
H^{k,p}=\sum_{n=k}^p \gamma_n h_n(x/R_n)
\end{equation}
and $\Gamma^{k,p}$ the skew-symmetric matrix given by $\Gamma^{k,p}_{1,2}(x_1,x_2)=H^{k,p}(x_1)$. Let $D(\Gamma^{0,p})$ be the effective diffusivity associated to homogenization of the periodic operator $L_{\Gamma^{0,p}}=1/2\Delta-\nabla \Gamma^{0,p} \nabla$. Then it is easy to see that
\begin{equation}\label{vjhsaassahv}
D(\Gamma^{0,p})=\begin{pmatrix}
1 & 0 \\
0 & D(\Gamma^{0,p})_{22}
\end{pmatrix}
\end{equation}
and it will be shown that
\begin{theorem}\label{IntSUMSFMSIHFDatha}
For
\begin{equation}\label{ajshvjdhvakhv1}
\epsilon = 4K_1 \big(\rho_{\min}(\gamma_{\min}-1)\big)^{-1} <1
\end{equation}
\begin{equation}\label{skjahbxkjbwxb}
1+4(1-\epsilon)\sum_{k=0}^{p}\gamma_k^2 \leq D(\Gamma^{0,p})_{22}\leq 1+4(1+\epsilon)\sum_{k=0}^{p}\gamma_k^2
\end{equation}
\end{theorem}
The super-diffusive behavior can be explained and controlled by a perpetual homogenization process taking place over the infinite number of scales $0,\ldots,n,\ldots$.
The idea of the proof of theorem \ref{IntSUMSFDThmsqdishgc1}  is to distinguish, when one tries to estimate \eref{hgdvsdckijeq71}, the smaller scales  which have already been homogenized ($0,\ldots,n_{ef}$ called effective scales), the bigger scales  which have not had a visible influence on the diffusion ($n_{dri},\ldots,\infty$ called drift scales because they will be replaced by a constant drift in the proof) and some intermediate scales that manifest their particular shapes in the behavior of the diffusion ($n_{ef}+1,\ldots,n_{dri}-1=n_{ef}+n_{per}$ called perturbation scales because they will enter in the proof as a perturbation of the homogenization process over the smaller scales).\\
The number of effective scales of is fixed by the mean squared displacement of $y_t.e_1$. Writing $n_{ef}(t)=\inf\{n: t \leq R_n^2\}$ one proves that
\begin{equation}\label{jhsasiuiuw798761}
\E[(y_t.e_2)^2]\sim D(\Gamma^{0,n_{ef}(t)})t
\end{equation}
Assume for instance $\R_n=\rho^n$ and $\gamma_n=\gamma^n$ then $n_{ef}(t)\sim \ln t/(2 \ln \rho)$ and $$\E[(y_t.e_2)^2]\sim t^{1+\frac{\ln \gamma}{\ln \rho}}$$
We remark that the quantitative control is sharper than the one associated to a perpetual homogenization on a gradient drift \cite{Ow00a}; this is explained by the fact that the number of perturbation scales is limited to only one scale with a divergence free drift. Nevertheless the main difficulty is to control the influence of this intermediate scale and the core of that control is based on the following mixing stochastic inequality (we write $\T^d_R=R \T^d$ the torus of dimension $d$ and side $R$)
\begin{proposition}\label{SFmsdle15jua1}
Let $R>0$ and $f,G\in (H^1(\T^1_R))^2$ such that $\int_0^1 f(y)dy=0$ and $\int_0^1 G(y)dy=0$, let $r>0$  and $t>0$
\begin{equation}
\Big|\E\big[G(b_t)\int_0^t \partial_1 f(rb_s)\,ds\big]\Big|\leq \|f\|_{L^2(\T^1_R)} \|G'\|_{L^2(\T^1_R)} 2 r^{-1}
\end{equation}
\end{proposition}

\subsection{Physical interpretation}
We want to emphasize that to obtain a super-diffusive behavior of the solution of \eref{SFsuperdiffstochdiffequ} of the shape $\E[y_t^2]\sim t^{1+\nu}$ with $\nu>0$ it is necessary to assume the exponential rate of growth of the parameters $\gamma_n$, this has a clear meaning when the flow is compared on a heuristic point of view to a real turbulent flow.\\
First, note that our model starts with the dissipation scale and expresses the inertial range of scales as geometrically
divergent. Next, the parameters $\gamma_n \|h_n'\|_\infty/R_n$ represents the amplitude of the pulsations of the eddies of size $R_n$.  It is a well known characteristic of  turbulence (\cite{Lali84} p. 129) that the amplitude of the pulsations increase with the scale, since for all scales $\|h_n'\|_\infty\leq  K_1$ one should have $\lim\gamma_n/R_n\rightarrow \infty $ to reflect that image. In our model, it is sufficient to assume the exponential rate of divergence of $\gamma_n$ to obtain a super-diffusive behavior.\\
 Let us also notice on a heuristic point of view (our model is not isotropic and does not depend on the time) that
the energy dissipated per unit time and unit volume in the eddies of scale $n$ is of order of
\begin{equation}
\epsilon_n \propto \frac{\gamma_n^2}{ R_n^4} K_2^2
\end{equation}
So saying that the energy is dissipated mainly in the small eddies is equivalent to saying that $\gamma_n/R_n^2 \rightarrow 0$ as $n \rightarrow \infty$ or if $R_n=\rho^n$ and $\gamma_n=\rho^{\alpha n}$, this equivalent to say that $\alpha<2$, which is included in the hypothesis \ref{csdcssjhxnjj9jdxc7h2}.\\
The Kolmogorov-Obukhov's law is equivalent to say that $K_2<\infty$, for all $n$ $h_n'(0)=0$ and
\begin{equation}\label{ksjdhdoihi8981}
\gamma_n \propto R_n^{\frac{4}{3}}
\end{equation}
or if $R_n=\rho^n$ and $\gamma_n=\rho^{\alpha n}$, this is equivalent to say that $\alpha=\frac{4}{3}$ corresponding to the Kolmogorov spectrum, which violates the hypothesis  \ref{csdcssjhxnjj9jdxc7h}, $\rho > \gamma$ but not the hypothesis \ref{csdcssjhxnjj9jdxc7h2} ($\rho>\gamma^\frac{1}{2}$) which will be addressed below.

\paragraph*{Overlapping ratios}
The super-diffusive behavior in theorem \ref{IntSUMSFDThmsqdishgc1} requires a minimal separation between scales, i.e. $\rho_{\min}>\rho_0$ and  this condition is necessary. Assume for instance $R_n=\rho^n$ and $\gamma_n=\gamma^n$,  then if $h_n(x_1)=g(x_1)-\gamma^{p} g(a^{p} x_1)$ (with $g\in C^1(\T^1)$ where $\T^1$ is the torus $\R/\Z$ and $a\in \N^*$)
it is easy to see that for $\rho=a$, $\Gamma$ is bounded and $y_t$ has a normal behavior by Norris's Aronson type estimates \cite{Nor97}. Thus, as for a gradient drift \cite{Ow00a}, it is easy to see on simple examples that when $\rho \leq \rho_0$ the solution of \eref{SFsuperdiffstochdiffequ} may have a normal or a super-diffusive behavior depending on the value of $\rho$ and the shapes of $h_n$ and ratios of normal behavior may be surrounded by ratios of anomalous behavior. This phenomenon is created by a strong overlap between spatial ratios between scales.\\
Thus the hypothesis \eref{akjsbdu8191} is necessary not only to ensure that one has a well defined $C^1$ stream matrix $\Gamma$ but also that its associated diffusion may show an anomalous fast behavior. Indeed with $h_n(x)=\sin(2\pi x)-3\sin(6 \pi x)$, $\gamma_n=\rho_n=3^n$ one has $\|h\|_\infty<\infty$, which leads to a normal behavior of $y_t$.\\
The case  $\rho_{\min}\leq \gamma_{\max}$ will be addressed with hypothesis \ref{csdcssjhxnjj9jdxc7h2}. Let us observe that to investigate on this case we had to add further informations on higher derivate of $h_n'$ to ensure that $h$ is well defined:
$h_n'(0)=0$, $\|h_n''\|_\infty\leq K_2$  under the constraint $\rho_{\min}> \gamma_{\max}^2$ (which includes the Kolmogorov case $\gamma=\rho^{4/3}$).

\subsection{Under the hypothesis \ref{csdcssjhxnjj9jdxc7h2}}
In this subsection we will give theorems putting into evidence the anomalous fast behavior of the solutions of \eref{SFsuperdiffstochdiffequ} under the hypothesis \ref{csdcssjhxnjj9jdxc7h2} which includes the Kolmogorov spectrum.

\begin{theorem}\label{djhsdkjdhb7710}
Under the hypotheses \eref{csdcssjhxnjj9jdxc7h2}, if there exists a constant $z>2$ such that $\gamma_{\min}\geq z$
then there exists a constant $C_0$ depending on $z,K_0,K_2$ such that if
\begin{equation}\label{djhsdkjdhb771}
\rho_{\min}^2 \gamma_{\max}^{-1} >C_0
\end{equation}
then there exist  constants $C_1>0$ depending on $z$, $C_2>0$ depending on $z,K_0,K_2,\gamma_{\max}$,  $C_3>0$ on $z,K_0$ and $C_4>0$ on $z,K_0,\rho_{\max}$ such that if $y_t$ is a solution of \eref{SFsuperdiffstochdiffequ} and $t\geq C_3$, then
\begin{equation}\label{djhsdkjdhb772}
C_1 t \gamma_{p(t)}^2 \leq \E\big[(y_t.e_2)^2\big] \leq C_2 t \gamma_{p(t)}^2
\end{equation}
and
\begin{equation}\label{djhsdkjdhb773}
C_4 t^{1+\nu(t)} \leq \E\big[(y_t.e_2)^2\big] \leq C_2 t^{1+\nu(t)}
\end{equation}
with
\begin{equation}\label{djhsdkjdhb774}
p(t):=\sup\{n \in \N\,:\, 16 (1+K_0^2) \big(1-(\gamma_{\min}-1)^{-1}\big)^{-2} R_n^2 <t\}
\end{equation}
and
\begin{equation}\label{djhsdkjdhb775}
\nu(t):=\ln \gamma_{p(t)}\big/ \ln R_{p(t)}
\end{equation}
\end{theorem}

\paragraph*{Remark}
It is important to note that the equation \eref{djhsdkjdhb772} emphasizes the role of never-ending homogenization in the fast behavior of the diffusion and its proof is also based on the separation of the infinite number of scales into smaller ones which act trough their effective properties, larger ones which will be bounded by a constant drift and an intermediate one that one has to control. The equation \eref{djhsdkjdhb773} gives a quantitative control of the anomaly without any self-similarity assumption; $\nu(t)$ is not a constant if the multi-scale velocity distribution associated to \eref{SFsuperdiffstochdiffequ} is not self-similar but one has
\begin{equation}\label{djhsdkjdhb77ss5}
\ln \gamma_{\min}\big/ \ln \rho_{\max}  \leq \nu(t)\leq \ln \gamma_{\max}\big/ \ln \rho_{\min}
\end{equation}

\begin{definition}\label{djhsdkjdhb7716}
The Stochastic Differential Equation \eref{SFsuperdiffstochdiffequ} is said to have a self-similar velocity distribution if and only if
\begin{equation}\label{djhsdkjdhb7717}
\rho_{\min}=\rho_{\max}=\rho
\end{equation}
and
\begin{equation}\label{djhsdkjdhb7718}
\gamma_{\min}=\gamma_{\max}=\gamma
\end{equation}
Then we will write
\begin{equation}\label{djhsdkjdhb7719}
\alpha=\ln \gamma/\ln \rho
\end{equation}
\end{definition}
From  theorem \ref{djhsdkjdhb7710} one easily deduces the following theorem.
\begin{theorem}\label{djhsdkjdhb77110}
Assume that the SDE \eref{SFsuperdiffstochdiffequ} has a self-similar velocity distribution with $0<\alpha<2$.
 Under the hypotheses \eref{csdcssjhxnjj9jdxc7h2},  there exists a constant $C_0$ depending on $\alpha,K_0,K_2$ such that if
\begin{equation}\label{djhsdkjdhb77111}
\rho >C_0
\end{equation}
then there exist  constants $C_1>0$ depending on $K_0,\rho$, $C_2>0$  on $K_0,K_2,\rho,\alpha$ and $C_3$ on $K_0$  such that if $y_t$ is a solution of \eref{SFsuperdiffstochdiffequ} and $t\geq C_3$, then
\begin{equation}\label{djhsdkjdhb77112}
C_1 t^{1+\alpha} \leq \E\big[(y_t.e_2)^2\big] \leq C_2 t^{1+\alpha}
\end{equation}
\end{theorem}
\paragraph*{Remark}
Observe that all the velocity distribution $0<\alpha<2$ is covered by theorem \ref{djhsdkjdhb77110}. The condition \eref{djhsdkjdhb77111} is still needed to avoid overlapping ratios and cocycles.\\
It is important to note that even if theorem \ref{djhsdkjdhb7710} allows to cover a larger spectrum of velocity distribution than theorem \ref{IntSUMSFDThmsqdishgc1}, its proof is based on the regularity of the drift associated to the diffusion ($K_2<\infty$). Whereas in theorem \ref{IntSUMSFDThmsqdishgc1}, the drift  $- \nabla \Gamma$ associated to the stochastic differential equation \eref{SFsuperdiffstochdiffequ} can be non Holder-continuous (only the regularity of $\Gamma$ is needed to prove an anomalous fast behavior)

\subsection{Remark: fast separation between scales}
When $\gamma_{\min}>1$ and $\gamma_{\max}<\infty$, the feature that distinguishes a strong anomalous behavior from a weak one is the rate at which spatial scales do separate. Indeed one can follow the proof of theorem \ref{IntSUMSFDThmsqdishgc1}, changing the condition $\rho_{\max}<\infty$ into  $R_n=R_{n-1}[\rho^{n^\alpha}/R_{n-1}]$ ($\rho,\alpha>1$)  and $\gamma_{\max}=\gamma_{\min}=\gamma$ to obtain
\begin{theorem}\label{IntSFSFDThmsqdishgc2}
For $t>t_0(\gamma_{2},R_2,K_0,K_1)$
\begin{equation}
C_1 t \gamma^{\beta(t)} \leq \E_0[|y_t.e_2|^2]\leq C_2 t \gamma^{\beta(t)}
\end{equation}
\begin{equation}
\text{with}\quad \quad\beta(t)= 2(2\ln \rho)^{-\frac{1}{\alpha}} (\ln t)^\frac{1}{\alpha}
\end{equation}
Where the constants $C_1$ and $C_2$ depends on $\rho,\gamma,\alpha,K_1$
\end{theorem}
And as $\alpha \downarrow 1$ the behavior of the solution of \eref{SFsuperdiffstochdiffequ} pass from weakly anomalous to strongly anomalous.

\section{Proofs under the hypothesis \ref{csdcssjhxnjj9jdxc7h}}\label{ksjhdbdshbh87h1}
\subsection{Notations and proof of theorem \ref{IntSUMSFMSIHFDatha}}
In this  subsection we will prove theorem \ref{IntSUMSFMSIHFDatha} using the explicit formula of the effective diffusivity. We will first introduce the basic notations that will also be used to prove theorem \ref{IntSUMSFDThmsqdishgc1}.
Let $J$ be smooth $T^2_1$ periodic  $2\times 2$ skew-symmetric matrix such that $J_{12}(x_1,x_2)=j(x_1)$ and consider
the periodic operator $L_J=1/2 \Delta-\nabla J\nabla$. We call $\chi_l$ the solution of the cell problem
associated  to $L_J$, i.e. the $T^2_1$ periodic solution of $L_J(\chi_l -l.x)=0$  with $\chi_l(0)=0$. One easily obtains that
\begin{equation}
\chi_l(x_1,x_2)=-2 l_2\Big[\int_0^{x_1} j(y) dy-x_1\int_0^1 j(y) dy\Big]
\end{equation}
The solution of the cell problem allows to compute the effective diffusivity $^tlD(J)l=\int_{T^2_1}|l-\nabla\chi_l(x)|^2dx$ that leads to
\begin{equation}\label{ggjaxgsvvvvvvvv}
D(J)=\begin{pmatrix}
1 & 0 \\
0 & 1+4\Var(j)
\end{pmatrix}
\end{equation}
For a $R$ periodic function $f$ we will write
\begin{equation}\label{dkdkjbbj998812}
\Var(f)=1/R \int_0^R (f(x)-\int_0^R f(y)dy)^2 dx
\end{equation}
Using the notation \eref{hvjshvxvajkko1}, from the equation \eref{ggjaxgsvvvvvvvv} we obtain \eref{vjhsaassahv} with
\begin{equation}\label{gdgcvjhgg12}
D(\Gamma^{0,p})_{2,2}=1+4\Var(H^{p})
\end{equation}
Now we will prove  by induction on $p$ that (for $\epsilon$ given by \eref{ajshvjdhvakhv1})
\begin{equation}\label{SFmatsupinHp6}
(1-\epsilon) \sum_{k=0}^p \gamma_k^2 \leq \Var(H^{p}) \leq (1+\epsilon) \sum_{k=0}^p \gamma_k^2
\end{equation}
Then the equation \eref{skjahbxkjbwxb} of theorem will follow by \eref{gdgcvjhgg12}.
The equation \eref{SFmatsupinHp6} is trivially true for $p=0$.\\
From the explicit formula of $H^{p}$ we will show in the paragraph \ref{shajhdvdksjhv71} that
\begin{equation}\label{SFmatsupinHp4}
\big|\Var(H^{p})-\Var(H^{p-1})-\gamma_p^2|\leq 2 \gamma_p K_1 r_p^{-1}\sqrt{\Var(H^{p-1})}
\end{equation}
Assuming \eref{SFmatsupinHp6} to be true at the rank $p$ one obtains
\begin{equation}\label{SFmatsupinHp455}
\begin{split}
\sqrt{\Var(H^{p})} &\leq (1+\epsilon)^\frac{1}{2} \gamma_{p+1} \Big(\sum_{k=0}^p (\gamma_k/\gamma_{p+1})^2 \Big)^\frac{1}{2}\\
&\leq (1+\epsilon)^\frac{1}{2}\gamma_{p+1}(\gamma_{\min}-1)^{-1}
\end{split}
\end{equation}
and it is easy to see that the condition \eref{ajshvjdhvakhv1} implies that
$\epsilon \geq 2K_1 (1+\epsilon)^\frac{1}{2} (\rho_{\min}(\gamma_{\min}-1))^{-1}$, combining this with \eref{SFmatsupinHp4} and \eref{SFmatsupinHp455} one obtains that
\begin{equation}
\big|\Var(H^{p+1})-\Var(H^{p})-\gamma_{p+1}^2| \leq \epsilon \gamma_{p+1}^2
\end{equation}
which proves the induction and henceforth the theorem.

\paragraph{}\label{shajhdvdksjhv71}
From the equation
$$\Var(H^{p})=\int_0^1 \big((H^{p-1}(R_p x)-\int_0^1 H^{p-1}(R_p y)dy)+\gamma_p (h_p(x)-\int_0^1 h_p(y) dy)\big)^2 dx
$$
one obtains
\begin{equation}\label{SFmatsupinHfgp4}
\big|\Var(H^{p})-\Var(H^{p-1})-\gamma_p^2|\leq 2 \gamma_p|E|
\end{equation}
with $$E= \Cov(S_{R_p}H^{p-1},h_p)$$ where $\Cov$ stands for the covariance: $\Cov(f,g)=\int_0^1 (f(x)-\int_0^1 f(y)dy)(g(x)-\int_0^1 g(y)dy)dx$ and $S_R$ is the scaling operator $S_R f(x)=f(Rx)$.\\
We will use the following mixing lemma whose proof is an easy exercise
\begin{lemma}\label{Pr_To_LD_separation2}
Let $(g,f)\in \big(C^1(T^{d}\big)^2$ and $R \in \N^*$ then
\begin{equation}
\Big |\int_{T^d_1}g(x) S_R f(x) dx -  \int_{T^d_1}g(x) dx \int_{T^d_1}f(x) dx\Big | \leq
(\|g'\|_\infty/R) \int_{T^d_1}\big|f\big| dx
\end{equation}
\end{lemma}
By lemma \ref{Pr_To_LD_separation2} and Cauchy-Schwartz inequality
\begin{equation}
\begin{split}\label{SFmatsupinHp3}
|E|\leq \frac{\|\nabla h_p\|_\infty R_{p-1}}{R_p}\int_0^1 |S_{R_p}H^{p-1}(x)|dx \leq \frac{K_1}{r_p}\sqrt{\Var(H^{p-1})}
\end{split}
\end{equation}
Combining this with \eref{SFmatsupinHfgp4} one obtains \eref{SFmatsupinHp4}.

\subsection{Anomalous mean squared displacement: theorem \ref{IntSUMSFDThmsqdishgc1}}
\subsubsection{Anomalous behavior from perpetual homogenization}\label{mchwlechlecj1}
Let $y_t$ be the solution of \eref{SFsuperdiffstochdiffequ}. Define
\begin{equation}\label{jhdsgkjdhg3122}
n_{ef}(t)=\inf\{n \in \N\,:\, t   \leq R_{n+1}^2 (\gamma_{n-1}/\gamma_{n+1})^2 2^{-3}K_1^{-2}(1-\gamma_{\max}/\rho_{\min})^2 \}
\end{equation}
$n_{ef}(t)$ will be the number of effective scales that have one can consider  \emph{homogenized} in the estimation of the mean squared displacement at the time $t$. Indeed, we will show in the sub subsection \ref{mchwlechlecj2} that
for $\rho_{\min}>C_{1,K_1,\gamma_{\min},\gamma_{\max}}$ and $t>C_{2,K_0,\gamma_{\min},\gamma_{1},R_1}$ one has
\begin{equation}\label{SFDmsdledegrtha1e22}
(1/4) t  \gamma_{n_{ef}(t)-1}^2 \leq \E\big[(y_t.e_2)^2\big]   \leq  C_{3,K_1,\gamma_{\min}} \gamma_{n_{ef}(t)+1}^2 t
\end{equation}
Combining this with the bounds \eref{Modsubboundrnrhonmin} and \eref{Modsubboundrnrhonmin002} on $R_n$ and $\gamma_n$ one obtains easily the theorem \ref{IntSUMSFDThmsqdishgc1}.

\subsubsection{Distinction between \emph{effective} and \emph{drift} scales. Proof of the equation \eref{SFDmsdledegrtha1e22}}\label{mchwlechlecj2}
By the Ito formula one has
\begin{equation}
y_t.e_2= \omega_t.e_2+\int_0^t \partial_1 h(\omega_s.e_1)ds
\end{equation}
And by  the independence of $\omega_t.e_2$ from $\omega_t.e_1$ one obtains
\begin{equation}\label{slakcxblkj81}
\E\big[(y_t.e_2)^2\big]= t+\E\Big[\big(\int_0^t \partial_1 h(\omega_s.e_1)ds\big)^2\Big]
\end{equation}
thus for all $p \in \N^*$, using $h=H^{p}+H^{p+1,\infty}$ one easily obtains (writing $H^p=H^{0,p}$)
\begin{equation}
\begin{split}
\E\big[(y_t.e_2)^2\big]&\geq  t+\frac{1}{2}\E\Big[\big(\int_0^t \partial_1 H^{p}(\omega_s.e_1)ds\big)^2\Big]-\E\Big[\big(\int_0^t \partial_1 H^{p+1,\infty}(\omega_s.e_1)ds\big)^2\Big]\\
&\leq  t+2\E\Big[\big(\int_0^t \partial_1 H^{p}(\omega_s.e_1)ds\big)^2\Big]+2\E\Big[\big(\int_0^t \partial_1 H^{p+1,\infty}(\omega_s.e_1)ds\big)^2\Big]\end{split}
\end{equation}
Now we will bound the larger scales $\partial_1 H^{p+1,\infty}$ as  drift scales, i.e. bound them by a constant drift using $\|\partial_1 H^{p+1,\infty}\|_\infty \leq K_1( \gamma_{p+1}/R_{p+1})(1-\gamma_{\max}/\rho_{\min})^{-1}$
\begin{equation}\label{jhsdghgsh61}
\begin{split}
\E\big[(y_t.e_2)^2\big]
&\geq t+\frac{1}{2}\E\Big[\big(\int_0^t \partial_1 H^{p}(\omega_s.e_1)ds\big)^2\Big]- \big(\frac{t K_1 \gamma_{p+1}}{R_{p+1}(1-\gamma_{\max}/\rho_{\min})})^2\\
&\leq t+2\E\Big[\big(\int_0^t \partial_1 H^{p}(\omega_s.e_1)ds\big)^2\Big]+2 \big(\frac{t K_1 \gamma_{p+1}}{R_{p+1}(1-\gamma_{\max}/\rho_{\min})})^2
\end{split}
\end{equation}
Write $I_p=\E\Big[\big(\int_0^t \partial_1 H^{p}(\omega_s.e_1)ds\big)^2\Big]$, since $H^{p}$ is periodic, for $t$ large enough, $I_p$ should behave like $t \gamma_p^2$. Nevertheless since the ratios between the scales are bounded, to control the asymptotic lower bound of the mean squared displacement in \eref{jhsdghgsh61}  we will need a quantitative control of $I_p$ that is sharp enough to show that the influence of the effective scales is not destroyed by the larger ones.
This control is based on stochastic mixing inequalities and
it will be shown in the sub subsection \ref{mchwlechlecj3} that for $\rho_{\min}> 8K_1/(\gamma_{\min}-1)$ one has
\begin{equation}\label{SFDmsdledegrtha1e22jhz781}
\begin{split}
I_p\leq &t 23  \big(\gamma_{p-1}^2 (1-1/\gamma_{\min})^{-2}+ K_0 K_1 \gamma_{p-1}(\gamma_p/r_p) (1-1/\gamma_{\min})^{-1}\big)
\\&+ t^2 K_1^2 (\gamma_p^2/R_p^2)
+ \sqrt{t} 68 \gamma_{p-1} \gamma_pR_{p-1}K_0^2(1-1/\gamma_{\min})^{-1}
\\&+ 60 K_0^2 (1-1/\gamma_{\min})^{-1} \gamma_{p-1}   \gamma_p (R_p^2 /r_p)
+16 K_0^2 \gamma_{p-1}^2 R_{p-1}^2
\end{split}
\end{equation}
and
\begin{equation}\label{SFDmsdledegrtha1e22jhz782}
\begin{split}
I_p\geq &t  \gamma_{p-1} (\gamma_{p-1}-16K_1 \gamma_p/r_p)
-\sqrt{t} 68 \gamma_{p-1}\gamma_pR_{p-1}K_0^2 (1-1/\gamma_{\min})^{-1}
\\&- 60 K_0^2 (1-1/\gamma_{\min})^{-1} \gamma_{p-1}   \gamma_p \frac{R_p^2}{r_p}
- 8 K_0^2 \gamma_{p-1}^2 R_{p-1}^2 (1-1/\gamma_{\min})^{-2}
\end{split}
\end{equation}
Choosing $p=n_{ef}(t)$ given by \eref{jhdsgkjdhg3122} one obtains  \eref{SFDmsdledegrtha1e22} from \eref{SFDmsdledegrtha1e22jhz781}, \eref{SFDmsdledegrtha1e22jhz782} and \eref{jhsdghgsh61} by straightforward computation under the assumption $\rho_{\min}>C_{K_0,K_1,\gamma_{\max},\gamma_{\min}}$.

\subsubsection{In{f}luence of the \emph{intermediate} scale on the \emph{effective} scales: proof of the inequalities \eref{SFDmsdledegrtha1e22jhz782} and \eref{SFDmsdledegrtha1e22jhz781}}\label{mchwlechlecj3}
In this sub subsection we will prove the inequalities \eref{SFDmsdledegrtha1e22jhz782} and \eref{SFDmsdledegrtha1e22jhz781} by distinguishing the scale $p$ as a perturbation scale, i.e. controlling its influence on the homogenization process over the
scales $0,\ldots,p-1$. More precisely,
writing $b_t$ the Brownian motion $\omega_t.e_1$ one has
\begin{equation}\label{sahxlhhloiu1}
I_p=I_{p-1}+\E\Big[\big(\int_0^t \partial_1 H^{p,p}(b_s)ds\big)^2\Big]+2J_p
\end{equation}
with
\begin{equation}\label{sahxlhhloiu2}
J_p=\E\Big[\big(\int_0^t \partial_1 H^{p,p}(b_s)ds\big)\big(\int_0^t \partial_1 H^{p-1}(b_s)ds\big)\Big]
\end{equation}
We will then control \eref{sahxlhhloiu1} by bounding $\partial_1 H^{p,p}$ by a constant drift to obtain
\begin{equation}\label{sahxlhhloiu3}
\E\Big[\big(\int_0^t \partial_1 H^{p,p}(b_s)ds\big)^2\Big]\leq t^2 K_1^2\gamma_p^2/R_p^2
\end{equation}
In the sub subsection \ref{mchwlechlecj4} we will control the homogenization process over the scales $0,\ldots,p-1$ and use our estimates on $D(\Gamma^{0,p-1})$ to obtain that for $\rho_{\min}> 8K_1/(\gamma_{\min}-1)$
\begin{equation}\label{sahxlhhloiu4}
I_{p-1}\leq t\gamma_{p-1}^2 20(1-1/\gamma_{\min})^2+ R_{p-1}^2\gamma_{p-1}^2 16K_0^2(1-1/\gamma_{\min})^2
\end{equation}
\begin{equation}\label{sahxlhhloiu5}
I_{p-1}\geq t\gamma_{p-1}^2 2(1-1/\gamma_{\min})^2- R_{p-1}^2\gamma_{p-1}^2 8K_0^2(1-1/\gamma_{\min})^2
\end{equation}
In the sub subsection \ref{mchwlechlecj5} we will bound $J_p$ is an error term by mixing stochastic inequalities to obtain
\begin{equation}\label{sahxlhhloiu6}
\begin{split}
|J_p|\leq \gamma_{p-1}\gamma_p(1-1/\gamma_{\min})^{-1} \Big( \sqrt{t}  R_{p-1}34 K_0^2   +(t/r_p) 8K_1  + R_p^2 r_p^{-1} 30 K_0^2 \Big)
\end{split}
\end{equation}
Combining \eref{sahxlhhloiu1}, \eref{sahxlhhloiu3}, \eref{sahxlhhloiu4}, \eref{sahxlhhloiu5} and \eref{sahxlhhloiu6} one obtains \eref{SFDmsdledegrtha1e22jhz781} and \eref{SFDmsdledegrtha1e22jhz782}

\subsubsection{Control of the homogenization process over the \emph{effective} scales: proof of the equations \eref{sahxlhhloiu4} and \eref{sahxlhhloiu5}}\label{mchwlechlecj4}
Writing for $m\leq p$,
\begin{equation}\label{saksjhvcxhcvh713}
\kappa^{m,p}=1/R_p\int_0^{R_p}H^{m,p}(y)dy
\end{equation}
 and $\kappa^p=\kappa^{0,p}$ one obtains by the Ito formula
\begin{equation}\label{meqchwlechlecj21}
\int_0^t \partial_1 H^{p-1}(b_s)ds=2\int_0^{b_t}(H^{p-1}(y)-\kappa^{p-1})dy-2\int_0^t (H^{p-1}(b_s)-\kappa^{p-1})db_s
\end{equation}
Using the periodicity of $H^p$:
\begin{equation}\label{meqchwlechlecj22}
\Big|\int_0^x \big(H^{p-1}(y)-\kappa^{p-1}\big)dy\Big|\leq R_{p-1} K_0 \gamma_{p-1} (1-1/\gamma_{\min})^{-1}
\end{equation}
Now we will show that for $\rho_{\min}> 8K_1/(\gamma_{\min}-1)$ and $p \in \N^*$
\begin{equation}\label{meqchwlechlecj23}
\begin{split}
\E[\int_0^t (H^{p}(b_s)-\kappa^{p})^2\,ds]&\leq (1-1/\gamma_{\min})^{-2} \big( t \gamma_p^2 5+ R_p^2 \gamma_p^2  4 K_0^2\big)\\ &\geq t \gamma_p^2 2 - R_p^2 \gamma_p^2  4 K_0^2  (1-1/\gamma_{\min})^{-2}
\end{split}
\end{equation}
Combining \eref{meqchwlechlecj21}, \eref{meqchwlechlecj22} and \eref{meqchwlechlecj23} one obtains \eref{sahxlhhloiu4} and \eref{sahxlhhloiu5} by straightforward computation.\\
The proof of \eref{meqchwlechlecj23} is based on standard homogenization theory: writing
\begin{equation}
f(x)= \int_0^x \big(H^p(y)-\kappa^p\big)^2- \Var(H^{p})x
\end{equation}
\begin{equation}
 \text{and}\quad g(x)=2 \Big(\int_0^x f(y)\,dy-(x/R_p)\int_0^{R_p}f(y)\,dy\Big)
\end{equation}
One obtains by Ito formula
\begin{equation}\label{dsjhdcdkhh1231}
\E[\int_0^t (H^{p}(b_s)-\kappa^{p})^2\,ds]=\Var(H^{p})t+\E[g(b_t)]
\end{equation}
Using the periodicity of $g$, one has $\|g\|_\infty \leq 4 K_0^2 \gamma_p^2 R_p^2 (1-1/\gamma_{\min})^{-2}$
Combing this with the estimate \eref{SFmatsupinHp6} on $\Var(H^p)$ one obtains \eref{meqchwlechlecj23} from \eref{dsjhdcdkhh1231}.

\subsubsection{Control of the in{f}luence of the  \emph{perturbation} scale: proof of the equation \eref{sahxlhhloiu6}}\label{mchwlechlecj5}
From the equations \eref{sahxlhhloiu2}, \eref{meqchwlechlecj21}
\begin{equation}
\Big|\int_0^x \big(H^{p-1}(y)-\kappa^{p-1}\big)dy\Big|\leq R_{p-1}\gamma_{p-1} K_0 (1-1/\gamma_{\min})^{-1}
\end{equation}
 and
\begin{equation}\label{meqchwlechlecj26541}
\int_0^t \partial_1 H^{p,p}(\omega_s.e_1)ds=2\int_0^{b_t}(H^{p,p}(y)-\kappa^{p,p})dy-2\int_0^t (H^{p,p}(b_s)-\kappa^{p,p})db_s
\end{equation}
One obtains using by straightforward computation  (using Cauchy-Schwartz inequality)
that
\begin{equation}\label{meqchwlecjlklecj2jk1}
\begin{split}
|J_p|\leq&  2\gamma_{p-1}\gamma_p R_{p-1}K_0^2 (1-1/\gamma_{\min})^{-1}  \sqrt{t}+4|J_{p,2}|+2|J_{p,3}|
\end{split}
\end{equation}
with
\begin{equation}
J_{p,2}=\E\Big[\int_0^t (H^{p-1}(b_s)-\kappa^{p-1})(H^{p,p}(b_s)-\kappa^{p,p})\,ds\Big]
\end{equation}
and
\begin{equation}
J_{p,3}=\E\Big[\int_0^t \partial_1 H^{p-1}(b_s)ds \int_0^{b_t} (H^{p,p}(y)-\kappa^{p,p})dy\Big]
\end{equation}
We will show in the sub subsection \ref{mchwlechlecj6} that the ratio $r_p$ allows a stochastic separation between the scales $0,\ldots,p-1$ and $p$ reflected by the following inequality
\begin{equation}\label{meqchwlecjlklecj2jk2}
\begin{split}
\big|J_{p,2}\big|\leq  \gamma_{p-1}\gamma_p K_0 (1-1/\gamma_{\min})^{-1}\big(8 K_0 R_{p-1}\sqrt{t}+2 t K_1/r_p\big)
\end{split}
\end{equation}
Using the proposition \ref{SFmsdle15jua1} (that we will prove in sub subsection \ref{mchwlechlecj7}) with
$G(x)=\int_0^{x} (H^{p,p}(y)-\kappa^{p,p})dy$, $r=r_p$ and $f(r_p x)=H^{p-1}(x)-k^{p-1}$ that
\begin{equation}\label{meqchwlecjlklecj2jk3}
\begin{split}
\Big|J_{p,3}\Big| \leq
 \gamma_{p-1}   \gamma_p R_p^2 r_p^{-1} 15 K_0^2 (1-1/\gamma_{\min})^{-1}
\end{split}
\end{equation}
Combining \eref{meqchwlecjlklecj2jk1}, \eref{meqchwlecjlklecj2jk2} and \eref{meqchwlecjlklecj2jk3} one obtains \eref{sahxlhhloiu6}.

\subsubsection{Stochastic separation between scales: proof of the equation \eref{meqchwlecjlklecj2jk2}}\label{mchwlechlecj6}
Writing for $x \in \R$,
\begin{equation}
g(x)=\int_0^x (H^{p-1}(y)-\kappa^{p-1})(H^{p,p}(y)-\kappa^{p}) dy
\end{equation}
one obtains by Ito formula
\begin{equation}\label{djhcvb71}
2\E[\int_0^{b_s} g(y)dy]=J_{p,2}
\end{equation}
Using the functional mixing lemma \ref{Pr_To_LD_separation2} one obtains easily that
\begin{equation}
\begin{split}
|g(x)| &\leq 2 K_0 \gamma_{p-1} (1-1/\gamma_{\min})^{-1}\big(4\gamma_p K_0 R_{p-1}+|x|K_1 (\gamma_p/r_p)\big)
\end{split}
\end{equation}
Combining this with \eref{djhcvb71} one obtains \eref{meqchwlecjlklecj2jk2}.

\subsubsection{Stochastic mixing: proof of the proposition \ref{SFmsdle15jua1}}\label{mchwlechlecj7}
By the scaling law of the Brownian motion  $$\E\big[G(b_t/R) \int_0^t (1/R)\partial_1 f(r b_s/R)\,ds\big]=R\E\big[G(b_{\frac{t}{R^2}})\int_0^{\frac{t}{R^2}} \partial_1 f(r b_s)\,ds\big]$$
and it is sufficient to prove the proposition assuming $R=1$. Let us write
\begin{equation}\label{jhgasgf}
I=\E\big[G(b_t)\int_0^t \partial_1 f(r b_s)\,ds\big]
\end{equation}
We will prove the proposition \ref{SFmsdle15jua1} by expanding \eref{jhgasgf} on the  the Fourier decompositions of $f$ and $G$ (written $f_k$ and $G_k$) and controlling trigonometric functions.
\begin{equation}
f(x)=\sum_{k\in \Z}f_k e^{i k 2\pi x}
\end{equation}
Write for $k,m \in \Z$
\begin{equation}
J_{k,m}=\int_0^t \E\big[e^{i k r 2\pi  b_s}e^{i m 2\pi  b_t}\big]\,ds
\end{equation}
By  straightforward computation
\begin{equation}
\begin{split}
J_{k,m}&=\int_0^t e^{-(2\pi)^2\frac{(kr+m)^2}{2}s-(2\pi)^2\frac{m^2}{2}(t-s)}\\
&= e^{-(2\pi)^2\frac{m^2}{2}t}\frac{1-e^{-(2\pi)^2(\frac{(kr)^2}{2}+krm)t}}{(2\pi)^2(\frac{(kr)^2}{2}+krm) }
\end{split}
\end{equation}
(in last fraction of the above equation, if the denominator is equal to $0$, we consider it as a limit to obtain the exact value $t$).\\
 Now
\begin{equation}\label{sahskjhdckb7861}
I=\sum_{k,m \in \Z^2} J_{km}ik 2\pi f_k G_m
\end{equation}
Thus since $f_0=G_0=0$, by Cauchy-Schwartz inequality
\begin{equation}\label{sahskjhdckb7862}
\begin{split}
|I|\leq \sum_{k \in \Z^*} |f_k| \big(\sum_{m \in \Z}(2\pi)^2m^2|G_m|^2\big)^\frac{1}{2}J_k^\frac{1}{2}
\end{split}
\end{equation}
with
\begin{equation}\label{sahskjhdckb7863}
J_k = \sum_{m \in \Z^*} \frac{1}{m^2} \Big(\frac{1-e^{-(2\pi)^2(\frac{(kr)^2}{2}+krm)t}}{(\frac{k^2r^2}{2}+krm)(2\pi)^2 }\Big)^2  e^{-(2\pi)^2 m^2 t}
\end{equation}
Using $(1-e^{-tx})/x \leq 3 t$ for $x>0$ and the fact that the minimum of $k^2r/2+km$ is reached for $m_0\sim kr/2$ and
we obtain
\begin{equation}\label{sahskjhdckb7864}
\begin{split}
J_k \leq& \sum_{m \in \Z^*}   e^{-(2\pi)^2 m^2 t^2} \big(\frac{4t }{k^2 r^2} + 2\sum_{m \in \Z^*} \frac{1}{m^2} \Big(\frac{1}{krm(2\pi)^2 }\Big)^2\big)\\
\leq & 4/r^2
\end{split}
\end{equation}
Observing that $\|G'\|_{L^2(\T^1)}^2= (2\pi)^2\sum_{m\in \Z}m^2 |G_m|^2$ one obtains from \eref{sahskjhdckb7862} and
\eref{sahskjhdckb7864}
\begin{equation}\label{sahskjhdckb7865}
\begin{split}
|I|\leq  \|G'\|_{L^2(\T^1)} (2/r)\sum_{k \in \Z^*} |f_k|
\end{split}
\end{equation}
 Which, using Cauchy-Schwartz inequality leads to
\begin{equation}\label{sahskjhdckb7866}
\begin{split}
|I|\leq  \|G'\|_{L^2(\T^1)}\|f\|_{L^2(\T^1)}  2/r
\end{split}
\end{equation}

\section{Proofs under the Hypothesis \ref{csdcssjhxnjj9jdxc7h2}}
In this section we will prove theorem \ref{djhsdkjdhb7710} under the hypothesis \ref{csdcssjhxnjj9jdxc7h2}. Observe that it is sufficient to prove the equation \eref{djhsdkjdhb772} under the hypotheses \ref{csdcssjhxnjj9jdxc7h2}, \eref{djhsdkjdhb771} and \eref{djhsdkjdhb774}.\\
We will use the same notations as in the section \ref{ksjhdbdshbh87h1}.
Let us observe that from the equation \ref{slakcxblkj81} one obtains
\begin{equation}\label{slakcxblkj8asx1}
\E\big[(y_t.e_2)^2\big]= t+\E\Big[X^2(b,0,t)\Big]+\E\Big[Y^2(b,0,t)\Big]+2\E\Big[X(b,0,t)Y(b,0,t)\Big]
\end{equation}
with for $p\in \N$ and $b_s=\omega_s.e_1$
\begin{equation}
X(b,0,t)=\int_0^t \partial_1 H^{0,p}(b_s)ds
\end{equation}
and
\begin{equation}
Y(b,0,t)=\int_0^t \partial_1 H^{p+1,\infty}(b_s)ds
\end{equation}
We will prove in the subsection \ref{kjdhsdbhhhhhh71} the following lemma (which is the core of the proof of theorem \ref{djhsdkjdhb7710}) which gives a quantitative control on decorrelation between the drifts associated with the small scales and those associated to the larger ones.
\begin{lemma}\label{sshbhh9j1j91}
For $t>R_p^2$
\begin{equation}\label{kdjdlskjhgjd070}
\begin{split}
\big|\E[X(b,0,t)Y(b,0,t)]\big|\leq &\big(t^{3/2}  R_{p+1}^2 R_{p+2}^{-2}+ 8 t R_p \big) R_{p+1}^{-2}R_p \gamma_p \gamma_{p+1}\\&  12 K_0 K_2 (1-\gamma_{\min}^{-1})^{-1} (1-\gamma_{\max}/\rho_{\min}^2)^{-1}
\end{split}
\end{equation}
\end{lemma}
Now by standard homogenization as it has been done in subsection \ref{mchwlechlecj4} it is easy to prove that
\begin{equation}\label{sajvsjhvdu81h1}
\E[X(b,0,t)^2]\leq  (R_p^2+t)\gamma_p^2 8 K_0^2   (1-\gamma_{\min}^{-1})^{-1}
\end{equation}
and
\begin{equation}\label{skjskjbsj891d1}
\E[X(b,0,t)^2]\geq  2 t \Var(H^p)- 8 R_p^2 \gamma_p^2 K_0^2   (1-\gamma_{\min}^{-1})^{-1}
\end{equation}
Now we will give in lemma \ref{jsdkdjdkjh81791} (proven in the subsection \ref{hvjslkk8ik1982}) the exponential rate of divergence of $\Var(H^p)$. \begin{lemma}\label{jsdkdjdkjh81791}
For $\gamma_{\min}>2$ one has
\begin{equation}\label{sskhksus8976}
\gamma_p^2 \big(1-(\gamma_{\min}-1)^{-1}\big)^2 \leq \Var(H^p)\leq \gamma_p^2 (1-\gamma_{\min}^{-1})^{-2}
\end{equation}
\end{lemma}
The proof of  lemma \ref{jsdkdjdkjh81791} is different from the one of theorem \ref{IntSUMSFMSIHFDatha} and is based on the domination of the influence of the biggest scale over the smaller ones (which is ensured by $\gamma_{\min}>2$, which is also necessary in the case $\rho_{\min}>\gamma_{\max}$ to avoid cocycles).\\
Combining \eref{sskhksus8976} with \eref{skjskjbsj891d1} one obtains that
\begin{equation}\label{skjskjbsj891d771}
\E[X(b,0,t)^2]\geq  2 t \gamma_p^2 \big(1-(\gamma_{\min}-1)^{-1}\big)^2 - 8 R_p^2 \gamma_p^2 K_0^2   (1-\gamma_{\min}^{-1})^{-1}
\end{equation}
Now let us observe that from \eref{slakcxblkj8asx1} and Cauchy-Schwartz inequality one obtains

\begin{equation}\label{slakcxblkj8ascccx1}
\begin{split}
\E\big[(y_t.e_2)^2\big]&\leq  t+2\E\Big[X^2(b,0,t)\Big]+2\E\Big[Y^2(b,0,t)\Big]\\
&\geq t+\E\Big[X^2(b,0,t)\Big]-2 \big|\E[X(b,0,t)Y(b,0,t)]\big|
\end{split}
\end{equation}
Combining \eref{skjskjbsj891d771}, \eref{sajvsjhvdu81h1}, \eref{kdjdlskjhgjd070}, \eref{slakcxblkj8ascccx1} and bounding $\E\Big[Y^2(b,0,t)\Big]$ by
\begin{equation}\label{slakcxblkj81221sasx1}
\begin{split}
\E\Big[Y^2(b,0,t)\Big]\leq& \|\partial_1^2 H^{p+1,\infty}\|_\infty^2 t^2 \E[b_t^2]\\
\leq & \gamma_{p+1}^2 t^3 R_{p+1}^{-4} \Big(K_2 (1-\gamma_{\min}^{-1})^{-1} (1-\gamma_{\max}/\rho_{\min}^2)^{-1}\Big)^2
\end{split}
\end{equation}
one obtains that for $t>R_p^2$, $\gamma_{\min}>2$ and $\rho_{\min}^2>\gamma_{\max}$ that
\begin{equation}\label{dshdsjhh891011}
\begin{split}
\E[(y_t.e_2)^2]\leq &32(R_p^2+t)\gamma_p^2  K_0^2   +8
\gamma_{p+1}^2 t^3 R_{p+1}^{-4} \Big(K_2  (1-\gamma_{\max}/\rho_{\min}^2)^{-1}\Big)^2
\end{split}
\end{equation}
and
\begin{equation}\label{dshdsjhh891012}
\begin{split}
\E[(y_t.e_2)^2]\geq & 2t \gamma_p^2 \big(1-(\gamma_{\min}-1)^{-1}\big)^2 - 16 R_p^2 \gamma_p^2 K_0^2
\\&-\big(t^{3/2}  R_{p+1}^2 R_{p+2}^{-2}+ 8 t R_p \big) R_{p+1}^{-2}R_p \gamma_p \gamma_{p+1} 48 K_0 K_2 (1-\gamma_{\max}/\rho_{\min}^2)^{-1}
\end{split}
\end{equation}
It follows \eref{dshdsjhh891012} and \eref{dshdsjhh891011} that for
\begin{equation}
16 (1+K_0^2) \big(1-(\gamma_{\min}-1)^{-1}\big)^{-2} R_p^2 <t \leq 16 (1+K_0^2) \big(1-(\gamma_{\min}-1)^{-1}\big)^{-2} R_{p+1}^2
\end{equation}
One has
\begin{equation}\label{ksjjskjndjjj891}
\begin{split}
\E[(y_t.e_2)^2]\geq &t \gamma_{p}^2 \Big(\big(1-(\gamma_{\min}-1)^{-1}\big)^{2}-\gamma_{\max}\rho_{\min}^{-2} 400 K_0K_2 (1-\gamma_{\max}/\rho_{\min}^2)^{-1}\\&\big(1+(1+K_0)\rho_{\min}^{-1}(\gamma_{\min}-1)(\gamma_{\min}-2)^{-1}\big)\Big)
\end{split}
\end{equation}
and
\begin{equation}\label{jhhjhdfhh97192}
\begin{split}
\E[(y_t.e_2)^2]\leq& t \gamma_{p}^2 (1+K_0^2)\\&\Big( 1+2050 \gamma_{\max}^2\big(1-(\gamma_{\min}-1)^{-1}\big)^{-4} K_2^2 (1-\gamma_{\max}/\rho_{\min}^2)^{-2}\Big)
\end{split}
\end{equation}
Which proves the equation \eref{djhsdkjdhb772} under the hypotheses \ref{csdcssjhxnjj9jdxc7h2}, \eref{djhsdkjdhb771} and \eref{djhsdkjdhb774} and thus theorem \ref{djhsdkjdhb7710}.

\subsection{Proof of  lemma \ref{sshbhh9j1j91}}\label{kjdhsdbhhhhhh71}
Let us introduce $Z$ a random variable idependant from $b_s$ and taking its values in $[0,R_p]$ uniformly with respect to the Lebesgue measure.
First, let us observe that by the periodicity of $H^p$
\begin{equation}\label{dsdkjndwjj89k981}
\E[X(b+Z,0,t)Y(b,0,t)]=0
\end{equation}
Next, define
$$\tau=\inf\{s>0\,:\, |b_s+Z|=0\;\text{or}\; |b_s+Z|=R_p\}$$
Now we will decompose $X(b+Z,0,t)$ and $Y(b,0,t)$ as follows,
\begin{equation}\label{skjhhbbsh}
X(b+Z,0,t)=X(b+Z,0,\tau)+ X(b+Z,\tau,t)
\end{equation}
\begin{equation}\label{djhbcldhbh2}
Y(b,0,t)=Y(b,0,\tau)+ Y(b,\tau,t)
\end{equation}
Let $b'$ be a BM independent from $b$ and $Z$ with the same law as $b$.
Combining the decompositions \eref{skjhhbbsh} and \eref{djhbcldhbh2} with the strong Markov property and the periodicity of $H^p$ one obtains that
\begin{equation}\label{kdjdlskjd7}
\begin{split}
\E[X(b+Z,0,t)Y(b,0,t)]=&\E[X(b+Z,0,\tau\wedge t) Y(b,0,\tau\wedge t)]\\&+\E[X(b',0,(t-\tau)_+)Y(b'+b_{\tau},0,(t-\tau)_+)]\\
&+\E[X(b+Z,0,\tau\wedge t)Y(b,\tau \wedge t,t)]\\&+ \E[X(b+Z,\tau\wedge t,t) Y(b,0,\tau\wedge t)]
\end{split}
\end{equation}
Where we have used $$\E[X(b',0,(t-\tau)_+)Y(b'+b_{\tau},0,(t-\tau)_+)]=\E[X(b+Z,\tau \wedge t,t)Y(b,\tau\wedge t,t)]$$
Using a similar decomposition for $\E[X(b,0,t)Y(b,0,t)]$ and subtracting \eref{kdjdlskjd7} combined with \eref{dsdkjndwjj89k981} one obtains that
\begin{equation}\label{kdjdlskjhgjd7}
\begin{split}
\E[X(b,0,t)Y(b,0,t)]=I_1+I_2+I_3+I_4+I_5+I_6+I_7
\end{split}
\end{equation}
With
\begin{equation}\label{sksjhblsh771}
I_1=\E\Big[X(b',0,(t-\tau)_+)\big(Y(b',0,(t-\tau)_+)-Y(b'+b_{\tau},0,(t-\tau)_+)\big)\Big]
\end{equation}
\begin{equation}\label{sksjhblsh772}
I_2=\E\Big[X(b',(t-\tau)_+,t)\big(Y(b',(t-\tau)_+,t)\Big]
\end{equation}
\begin{equation}\label{sksjhblsh773}
I_3=-\E\Big[X(b+Z,0,\tau\wedge t) Y(b,0,\tau\wedge t)\Big]
\end{equation}
\begin{equation}\label{sksjhblsh774}
I_4=-\E\Big[X(b+Z,0,\tau\wedge t)Y(b,\tau \wedge t,t)\Big]
\end{equation}
\begin{equation}\label{sksjhblsh775}
I_5=-\E\Big[X(b+Z,\tau\wedge t,t) Y(b,0,\tau\wedge t)\Big]
\end{equation}
\begin{equation}\label{sksjhblsh776}
I_6=\E\Big[X(b',0,(t-\tau)_+) Y(b',(t-\tau)_+,t)\Big]
\end{equation}
\begin{equation}\label{sksjhblsh777}
I_7=\E\Big[Y(b',0,(t-\tau)_+) X(b',(t-\tau)_+,t)\Big]
\end{equation}
We shall use standard homogenization techniques to estimate the influence of the effective scales in those terms.
Let us observe that by the Ito formula one has a.s. for all $t>0$,
\begin{equation}\label{shhjhhhhi71}
X(b,0,t)=\chi^p(b_t)-2\int_0^t \big(H^p(b_s)-\kappa^p\big)db_s
\end{equation}
with (using the notation introduced in \eref{saksjhvcxhcvh713})
\begin{equation}\label{shhjhhhhi7102}
\chi^p(x)=2\int_0^{x} H^p(y)dy-x 2\kappa^p
\end{equation}
Next we will control the larger scales by bounding the growth their drift. Indeed
Using the  uniform control on the second derivate of the functions $h_n$ in the hypothesis \ref{csdcssjhxnjj9jdxc7h2} one has
\begin{equation}\label{hhdkkks89891}
\|\partial_1^2 H^{p+1,\infty}\|_\infty\leq K_2 \gamma_{p+1}R_{p+1}^{-2} (1-\gamma_{\max}/\rho_{\min}^2)^{-1}
\end{equation}
Using Cauchy-Schwartz inequality in \eref{sksjhblsh772}, \eref{shhjhhhhi71} and \eref{hhdkkks89891} one obtains that
\[
\begin{split}
|I_2|\leq & \Big(\E\big[X(b',0,(t-\tau)_+)^2\big]\Big)^\frac{1}{2}\Big(\E\big[(Y(b',(t-\tau)_+,t))^2\big]\Big)^\frac{1}{2}\\
\leq & (1-\gamma_{\min}^{-1})^{-1}\big( 4 K_0^2 \gamma_p^2 R_p^2+4 K_0^2 \gamma_p^2 t\big)^\frac{1}{2}
\|\partial_1^2 H^{p+1,\infty}\|_\infty \Big(\E\big[|b'_t|^2 \tau^2\big]\Big)^\frac{1}{2}
\end{split}
\]
Thus for
\begin{equation}
t>R_p^2
\end{equation}
\begin{equation}\label{skjhddbobu8772}
\begin{split}
|I_2| \leq  \sqrt{t}  \gamma_p \gamma_{p+1} R_p^3 R_{p+1}^{-2} 12 K_0 K_2 (1-\gamma_{\min}^{-1})^{-1} (1-\gamma_{\max}/\rho_{\min}^2)^{-1}
\end{split}
\end{equation}
Similarly using \eref{shhjhhhhi71} one obtains that
\begin{equation}\label{skjhddbobu8773}
\begin{split}
|I_3|\leq   \gamma_p \gamma_{p+1} R_p^4 R_{p+1}^{-2} 12 K_0 K_2(1-\gamma_{\min}^{-1})^{-1}  (1-\gamma_{\max}/\rho_{\min}^2)^{-1}
\end{split}
\end{equation}
To estimate $I_4$ we use  the conditional independence
\begin{equation}\label{skjslkdjd761}
\begin{split}
I_4=&-\E\Big[1_{\tau\wedge t<t}\E\big[X(b+Z,0,\tau\wedge t)\big|\tau \wedge t\big]\E\big[Y(b,\tau \wedge t,t)\big|\tau \wedge t\big]\Big]\\
&-\E\Big[1_{\tau\wedge t=t}\E\big[X(b+Z,0,\tau\wedge t)\big|\tau \wedge t\big]\E\big[Y(b,\tau \wedge t,t)\big|\tau \wedge t\big]\Big]
\end{split}
\end{equation}
Since $Y(b, t,t)=0$ the second term in \eref{skjslkdjd761} is null. Moreover conditionally on the event $\tau <t$, $b_{\tau\wedge t}+Z$ with equal chances $1/2$ is equal to $0$ or $R_p$, so using \eref{shhjhhhhi71} one obtains
that on $\{\tau\wedge t<t\}$, $\E\big[X(b+Z,0,\tau\wedge t)\big|\tau \wedge t\big]=0$, which leads to
\begin{equation}\label{skjhddbobu8774}
\begin{split}
I_4=0
\end{split}
\end{equation}
Similarly, using conditional independence \eref{shhjhhhhi71} and \eref{hhdkkks89891} one easily obtains
\begin{equation}\label{skjhddbobu8775}
\begin{split}
|I_5|\leq \gamma_p \gamma_{p+1} R_p^4 R_{p+1}^{-2} 12 K_0 K_2  (1-\gamma_{\min}^{-1})^{-1}(1-\gamma_{\max}/\rho_{\min}^2)^{-1}
\end{split}
\end{equation}
\begin{equation}\label{skjhddbobu8776}
\begin{split}
|I_6|\leq \sqrt{t} \gamma_p \gamma_{p+1} R_p^3 R_{p+1}^{-2} 12 K_0 K_2 (1-\gamma_{\min}^{-1})^{-1} (1-\gamma_{\max}/\rho_{\min}^2)^{-1}
\end{split}
\end{equation}
To estimate $I_7$ let us observe that from \eref{shhjhhhhi71} and the spectral gap associated to the Laplacian on the torus (which is equivalent to Poincare inequality, see \cite{Ol94}) one has
\begin{equation}
\begin{split}
|\E[X(b',(t-\tau)_+,t)|t-\tau]|&\leq |\chi^p(b_t)-\chi^p(b_{(t-\tau)_+})|
\\&\leq 12 e^{-R_p^{-2}/2 (t-\tau)_+} K_0 (1-\gamma_{\min}^{-1})^{-1}\gamma_p R_p
\end{split}
\end{equation}
Thus using conditional independence and distinguishing the events $\tau>t/2$ (whose probability is bounded using the spectral gap) and $\tau \leq t/2$ , one obtains
\begin{equation}\label{skjhddbobu8777}
\begin{split}
|I_7|&\leq t^{3/2} e^{-t R_p^{-2}/4} \gamma_p \gamma_{p+1}16 R_p R_p^{-2} K_0 K_2 (1-\gamma_{\min}^{-1})^{-1} (1-\gamma_{\max}/\rho_{\min}^2)^{-1}
\\& \leq \gamma_p \gamma_{p+1}16 R_p^4 R_p^{-2} K_0 K_2  (1-\gamma_{\min}^{-1})^{-1}(1-\gamma_{\max}/\rho_{\min}^2)^{-1}
\end{split}
\end{equation}
To estimate $I_1$ we will distinguish the scale $p+1$ as an intermediate scale. First observe that
\begin{equation}\label{sksjhblsh7711101}
Y(b',0,(t-\tau)_+)-Y(b'+b_{\tau},0,(t-\tau)_+)=2\gamma_{p+1} J_1+2J_2
\end{equation}
with
\begin{equation}
J_2=\int_0^{(t-\tau)_+} \big(\partial_1 H^{p+2,\infty}(b_s')-\partial_1 H^{p+2,\infty}(b_s'+b_\tau)\big)ds
\end{equation}
and
\begin{equation}
J_1=\int_0^{(t-\tau)_+} \big(\partial_1 h_{p+1}(b_s')-\partial_1 h_{p+1}(b_s'+b_\tau)\big)ds
\end{equation}
Using $|b_\tau|\leq R_p$ and \eref{hhdkkks89891}, one obtains that for $t>R_p$
\begin{equation}\label{kedhxiiu771}
|\E\Big[X(b',0,(t-\tau)_+)J_2\Big]|\leq t^{3/2} R_p \gamma_p \gamma_{p+1} R_{p+2}^{-2} 2 K_0 K_2  (1-\gamma_{\min}^{-1})^{-1}(1-\gamma_{\max}/\rho_{\min}^2)^{-1}
\end{equation}
Now using Ito formula one has
\begin{equation}
\begin{split}
J_1&=\int_0^{b_{\tau}} h_{p+1}(x)dx-\int_{b_{(t-\tau)_+}'}^{b_{(t-\tau)_+}'+b_\tau} h_{p+1}(x)dx\\&-\int_0^{(t-\tau)_+}
\big(h_{p+1}(b_s')-h_{p+1}(b_s'+b_\tau)\big)db_s'
\end{split}
\end{equation}
Thus Observing that $(\tau,b_\tau)$ has the same law as $(\tau,-b_\tau)$ one obtains that
\begin{equation}\label{hjdhvwvi661}
\begin{split}
\E\Big[X(b',0,(t-\tau)_+)J_1\Big]=&\E\Big[X(b',0,(t-\tau)_+)\int_0^{b_{\tau}} \big(h_{p+1}(x)-h_{p+1}(-x)\big)/2\, dx\Big]\\
+\E\Big[X(b',0,(t-\tau)_+)&\int_{b_{(t-\tau)_+}'}^{b_{(t-\tau)_+}'+b_\tau} \big(h_{p+1}(x)-h_{p+1}(x-b_\tau)\big)/2 \,dx\Big]
\\-\E\Big[X(b',0,(t-\tau)_+)&\int_0^{(t-\tau)_+}
\big(h_{p+1}(b_s')-h_{p+1}(b_s'+b_\tau)\big)db_s'\Big]
\end{split}
\end{equation}
It follows from \eref{hjdhvwvi661} that for $t>R_p$
\begin{equation}\label{hjdhvwvi661222}
\begin{split}
\big|\E\Big[X(b',0,(t-\tau)_+)\gamma_{p+1}J_1\Big]\big|\leq 6 K_0 \gamma_p t R_p^2 \gamma_{p+1}R_{p+1}^{-2}  K_2 (1-\gamma_{\min}^{-1})^{-1} (1-\gamma_{\max}/\rho_{\min}^2)^{-1}
\end{split}
\end{equation}
Thus from \eref{hjdhvwvi661222} and \eref{kedhxiiu771} we deduce that
\begin{equation}\label{skjhddbobu8771}
|I_1|\leq \big(t^{3/2}  R_{p+2}^{-2}+ t R_p R_{p+1}^{-2}\big)R_p \gamma_p \gamma_{p+1}  12 K_0 K_2  (1-\gamma_{\min}^{-1})^{-1}(1-\gamma_{\max}/\rho_{\min}^2)^{-1}
\end{equation}
Combining \eref{skjhddbobu8771},\eref{skjhddbobu8772},\eref{skjhddbobu8773},\eref{skjhddbobu8774},\eref{skjhddbobu8775},\eref{skjhddbobu8776},\eref{skjhddbobu8777} and \eref{kdjdlskjhgjd070} we have obtained that for $t>R_p^2$ the equation \eref{kdjdlskjhgjd070} is valid.

\subsection{Proof of  lemma \ref{jsdkdjdkjh81791}}\label{hvjslkk8ik1982}
Observe that (using the notation introduced in \eref{saksjhvcxhcvh713})
\begin{equation}
\Var(H^p)=\Var(H^{p-1})+\gamma_p^2+2  R_p^{-1}\int_0^{R_p}\big(H^{p-1}(x)-\kappa^{0,p-1}\big)\big(h_p(x)-\kappa^{p,p}\big)
\end{equation}
It follows by Cauchy Schwartz inequality and the normalization condition \eref{SFsupallngamnzeroo} that
\[
\Var(H^{p-1})+\gamma_p^2+2 \gamma_p \Var(H^{p-1})^\frac{1}{2}\geq \Var(H^{p})\geq \Var(H^{p-1})+\gamma_p^2-2 \gamma_p \Var(H^{p-1})^\frac{1}{2}
\]
From which we deduce
\begin{equation}\label{skjsbbjj81}
 \Var(H^{p})^\frac{1}{2} \leq \gamma_p+\Var(H^{p-1})^\frac{1}{2}
\end{equation}
and
\begin{equation}\label{skjsbbjj82}
 \Var(H^{p})^\frac{1}{2} \geq \gamma_p-\Var(H^{p-1})^\frac{1}{2}
\end{equation}
From \eref{skjsbbjj81} one obtains by induction that
\begin{equation}\label{skjsbbjj83}
\Var(H^p)^\frac{1}{2} \leq \gamma_p (1-\gamma_{\min}^{-1})^{-1}
\end{equation}
Combining \eref{skjsbbjj82} with \eref{skjsbbjj83} one obtains  that
\begin{equation}\label{skjsbbjj84}
\Var(H^p)^\frac{1}{2} \geq \gamma_p \big(1-(\gamma_{\min}-1)^{-1}\big)
\end{equation}
And observe that $1-(\gamma_{\min}-1)^{-1}>0$ for $\gamma_{\min}>2$

\paragraph*{Acknowledgments}
The authors would like to thank the referees for useful comments. Part of this work was supported by the Aly Kaufman fellowship.

\end{document}